\documentclass[10pt]{article}
\usepackage[utf8]{inputenc}
\usepackage{latexsym,amssymb,amsmath,url,authblk}

\def\N{\mathbb{N}}

\def\Z{\mathbb{Z}}
\def\R{\mathbb{R}}

\def\proof{\par\noindent{\em Proof. }}
\def\eproof{\hfill{$\Box$}\bigskip}

\def\ds{\dots}
\def\sus{\subset}

\newtheorem{thm}{Theorem}[section]
\newtheorem{propo}[thm]{Proposition}
\newtheorem{corollary}[thm]{Corollary}
\newtheorem{lemma}[thm]{Lemma}
\newtheorem{prob}[thm]{Problem}
\newtheorem{defi}[thm]{Definition}

\begin{document}

\title{Rethinking real numbers as infinite decimals}
\markright{Decimal real numbers}
\author{Martin Klazar}
\affil{{\it Department of Applied Mathematics, Faculty of Mathematics and Physics, Charles University, Malostransk\'e n\'am\v est\'\i\ 25, 118 00 Praha 1, Czechia}\\ e-mail: {\tt klazar@kam.mff.cuni.cz}}

\maketitle

\begin{abstract}
We give a~detailed construction of the complete ordered field of real numbers by means of infinite decimal expansions. 
We prove that in the canonical encoding of decimals neither addition 
nor multiplication is {\em computable}, but that both operations are {\em weakly computable}; we introduce both kinds of computability in greater generality. 
We determine which additive and multiplicative shifts (restrictions of addition and multiplication to one variable) are computable, and 
prove that each of these shifts
becomes computable after a~permutation of encoding. We ask if it is the case for the bivariate addition and multiplication.
\end{abstract}

\section{Introduction. }In \cite{klaz} we discuss Dedekind's and, tangentially, Cantor's constructions of the real numbers, and unify them by means of
so called generating sequences in a~single construction. In the present companion article we apply this technique
to the third main construction of $\mathbb{R}$, usually attributed to {\em Karl Weierstrass (1815--1897)}\,---\,see 
\cite{encyR} or \cite{ljas_al}, based on infinite decimal expansions. Real numbers as infinite decimals are discussed, for example, 
in \cite{abia}, \cite{fard_li}, \cite{gowe_blog} and \cite[Chapter~4.1]{gowe}, \cite{ljas_al} and \cite{stol_gmei}. K.~Weierstrass lectured on 
real\,---\,in the contemporaneous terminology ``irrational''\,---\,numbers in his course on analytic functions in Berlin first in 1863--64
(\cite[p.~49]{twed}). Accounts of his construction were given for example in \cite{dant}, and recently in \cite{twed}.  
The history of old constructions of $\mathbb{R}$ was reviewed in \cite[pp. 396--400]{cajo}, and \cite{weis} surveys
in technical terms both old and recent constructions.

In our article we construct elements of $\R$ as infinite decimal expansion. They have an even more complicated arithmetic 
than Dedekind's cuts, and we improve it by the same device as in \cite{klaz}, by generating
sequences. We introduce their version for decimals below after the proof of Proposition~\ref{zeroth_res}. 
Our idea is to add and multiply decimals not by using their digits, but by means of their definitions in terms 
of Cauchy sequences of fractions. The satisfaction of axioms of an ordered field in $\mathbb{R}$ is then immediate because 
it transfers directly from the ordered field $\mathbb{Q}$.

The decimal construction of $\R$ is important because it produces real numbers in the form in which
they are used in practice. It is not so well known that one can view Cantor's construction as disqualified by the fact that 
it uses uncountable sets for representations of individual real numbers (\cite[Section~4]{klaz}). 
Dedekind's cuts are countable and even hereditarily at most countable (\cite[Section~4]{klaz}), 
but in practice one does not compute with them but with decimal expansions. Thus although the construction of $\R$ via decimals
is the last one in theoretical development of the three constructions, in terms of practice it is the most important one. 

We briefly survey content of our article. 
Section~2 contains definitions related to decimals. In Propositions~\ref{first_res}--\ref{third_res} we prove three key results 
on generating sequences of decimals, analogous to the
results in \cite[Section~2]{klaz} for cuts. Before that we show in Proposition~\ref{supremum} 
that two linear orders, that of decimals and that of decimals unified with false decimals, have suprema and infima,
and in Proposition~\ref{zeroth_res} that each bounded monotone sequence of terminating decimals converges to a~decimal.
These results are needed for the definition of generating sequences. 

In Section~3 we construct, by means of decimals, the real numbers as a~complete ordered field. The construction itself as given in the proof 
of Theorem~\ref{dede} is straightforward, but technical preparations for it in Section~2 and in Lemma~\ref{division} make it on the whole
more complicated than the construction by Dedekind's cuts in \cite{klaz}.

Sections~4 and 5 contain the main results of our article. They concern computing the addition $+_R$ and the multiplication $\cdot_R$ of decimals. This is
an attractive but also tantalizing topic that is difficult to grasp, we made already one attempt \cite{klaz2009} in 2009. In Section~4 we introduce in Definition~\ref{defi_comput} computability 
of $+_R$ and $\cdot_R$
in general, for any binary operation $$F\colon X\times X\to X$$ on a~(typically uncountable) set $X$ of infinite words over a finite alphabet. 
This setting includes besides the operations $+_R$ and $\cdot_R$ also, for example, the arithmetic of $p$-adic numbers: in Theorem~\ref{thm_comp} we show
that  in the defined sense their addition and multiplication is indeed computable. In contrast, by Theorem~\ref{thm_incomp} neither of the operations $+_R$ and $\cdot_R$ on decimals can be computed.
Corollary~\ref{cor_gen_coding} shows that this also holds in any encoding of decimals that is computably equivalent to the canonical one. It is
well known, though, that $+_R$ and $\cdot_R$ are computable in some weaker, intuitive sense. We formalize it in Definition~\ref{defi_weak_comp} in the notion of weak computation, and confirm the intuition in 
Theorems~\ref{thm_weak_comp_addi} and \ref{thm_weak_comp_mult}. 

It is easier to investigate addition and multiplication of decimals when they are restricted to one variable, and
this is the topic of Section~5. In Propositions~\ref{prop_on_add_shifts} and \ref{prop_on_mult_shifts} we determine which of the restrictions are computable, 
and in Theorems~\ref{thm_add_shifts} and \ref{thm_mul_shifts} we show that each restriction can be made computable by permuting the underlying encoding.

The last Section~6 poses two questions. 1. Which graphs of unary operations on $X$ are realized by computable operations? 2. Is it possible
to rename the underlying encoding so that the binary operation $+_R$, or $\cdot_R$, on decimals becomes computable?

\section{Decimals}

Let $\mathbb{Z}=\{\dots,-1,0,1,\dots\}$ be the integers, $\mathbb{N}=\{1,2,3,\dots\}$ be the natural numbers, $\mathbb{N}_0=\{0,1,2,\dots\}$ 
be the nonnegative 
integers, $\mathbb{Q}$ be the rational numbers (fractions), and $D$ be the ten-element set $D=\{0,1,\dots,9\}\subset\mathbb{N}_0$. 
A~{\em decimal $d$} is a~formal infinite linear combination
\begin{eqnarray*}
d&=&(-)\sum_{n\le k}d_n\cdot 10^n=(-)\sum_{n\le k}d_n 10^n=(-)(d_k10^k+_D d_{k-1}10^{k-1}+_D\cdots)\;,\\
&&\text{where $k\in\mathbb{N}_0$, $n\in\mathbb{Z}$ and $d_n\in D$}\;,
\end{eqnarray*}
that has the optional minus sign $-$ (the plus sign $+$ is not written) and satisfies three conditions: (i) if $d_k=0$ then $k=0$, (ii) never $d_n=d_{n-1}=\dots=9$ for any 
$n\le k$, and (iii) the decimal $d=-\sum_{n\le0}0\cdot10^n$ is not allowed. The symbol $+_D$ denotes formal 
addition; the true addition in $\mathbb{Q}$ is denoted by $+$.
We call the number $k$ the {\em order of $d$} and call each number $d_n$
the {\em digit at $10^n$ (in $d$)}. If need arises to refer to digits $d_n$ with $n>k$, we set them to zero. Two decimals are equal if and only 
if they have equal signs, equal orders, and for every $n\in\mathbb{Z}$ have equal digits at $10^n$.

We denote the set of decimals by $\mathbb{R}$. A~{\em terminating decimal} is one that from some $n$ below has zero digits, $d_n=d_{n-1}=\dots=0$. 
We denote the set of terminating decimals by ${\mathbb{R}}_{\mathrm{t}}\subset\mathbb{R}$. The {\em zero decimal} is 
$$
0_R:=\sum_{n\le0}0\cdot10^n\in\mathbb{R}_{\mathrm{t}}\;. 
$$
If $d\in\mathbb{R}_{\mathrm{t}}$ is nonzero and $n\in\mathbb{Z}$ is maximum with $d_{n-1}=d_{n-2}=\dots=0$, we write $d$ as a~finite formal linear combination
$d=(-)\sum_{i=n}^k d_i10^i$ (so $d_n\ne0$). Instead of $1\cdot 10^n$ we may write just $10^n$.
A~decimal $d$ is {\em positive} (resp. {\em negative}) if it is different from $0_R$ and does not have the minus sign (resp. has the minus sign). 
A~decimal $d$ is {\em nonnegative} if it is $0_R$ or is positive.
If $d$ is a~nonzero decimal, we obtain from it the decimal $d_-$ simply by omitting the minus sign if $d$ is negative, and adding it to $d$ if $d$ is positive. 
We set $(0_R)_-:=0_R$,

The {\em false decimals} form the set $\mathbb{R}_{\mathrm{f}}$, disjoint to $\mathbb{R}$, and arise when the above conditions (ii) and (iii) are
omitted. Thus 
$$
\mathbb{R}_{\mathrm{f}}=\{\overline{d}\;|\;d\in\mathbb{R}_{\mathrm{t}}\}\;,
$$
where $\overline{0_R}:=-\sum_{n\le0}0\cdot10^n$ and if $d=(-)\sum_{i=n}^k d_i10^i\ne0_R$ 
is a~terminating decimal (so $n\le k$ and $d_n\ne0$) then 
$$
\overline{d}:=(-)\sum_{m\le k'} e_m10^m\;,
$$
where $\overline{d}$
has the same sign as $d$, has order $k'=k$ or $k'=k-1$, and has the digits $e_m=9$ for every $m<n$, $e_n=d_n-1$, and $e_m=d_m$ for $m>n$.
If $n=k>0$ and $d_n=1$ then $k'=k-1$ (to keep condition (i)), else $k'=k$. The map $\overline{\cdots}\colon\mathbb{R}_{\mathrm{t}}\to\mathbb{R}_{\mathrm{f}}$ is 
a~bijection, and we denote its inverse also by $\overline{\cdots}$. False decimals arise as formal limits of (true) decimals. 
We extend in the obvious way to $\mathbb{R}\cup\mathbb{R}_{\mathrm{f}}$ the terminology of positive, negative and nonnegative (false) decimals, and also extend to $\mathbb{R}\cup\mathbb{R}_{\mathrm{f}}$ the map $d\mapsto d_-$, which we change for (false) zero decimals to $(0_R)_-:=-\sum_{n\le0}0\cdot10^n$ and $(-\sum_{n\le0}0\cdot10^n)_-:=0_R$.

We assign to each $d\in\mathbb{R}_{\mathrm{t}}$ the rational number $r(d)\in\mathbb{Q}$, obtained by evaluating $d$ in the usual arithmetic of $\mathbb{Q}$; we interpret the symbols in $d$ as elements of $\mathbb{Q}$ in the obvious way.
For example, 
$${\textstyle
r(-(2\cdot10^1+_D0\cdot10^0+_D3\cdot10^{-1}))=-20.3=-\frac{203}{10}\;.
}
$$
The map $r$ is injective and the subset 
$$
\mathbb{Q}_D:=r(\mathbb{R}_{\mathrm{t}})\subset\mathbb{Q}
$$
is a~subring of the field $\mathbb{Q}$. However, $\mathbb{Q}_D$ is in general not closed to division. For every $d,e\in\mathbb{R}_{\mathrm{t}}$ we define the binary operation
$+_R\colon\mathbb{R}_{\mathrm{t}}\times\mathbb{R}_{\mathrm{t}}\to\mathbb{R}_{\mathrm{t}}$ by 
$$
d+_R e:=r^{-1}(r(d)+r(e))\;,
$$
and define similarly $-_R$.
For any decimal $d=(-)\sum_{n\le k}d_n10^n$ and any $m\in\mathbb{N}$ we denote by 
$$
d\,|\,m\in\mathbb{R}_{\mathrm{t}}
$$ 
the {\em $m$-th truncation of $d$}, the terminating decimal $(-)\sum_{n\le k}e_n10^n$ that has the same sign and order as $d$, and has digits $e_n=d_n$ for $n\ge-m$, and $e_n=0$ for $n<-m$.

Recall that a~relation $<_X\,\sus X\times X$ is a~linear order on a~set $X$ if it is transitive 
($a<_X b$ and $b<_X c$ imply $a<_X c$), irreflexive (never $a<_X a$) and total (for any $a,b\in X$ either $a<_X b$ or $a=b$ or $b<_X a$). The notation $a\le_X b$
means that $a<_X b$ or $a=b$. We denote the standard linear order on $\mathbb{Q}$ by $<$. We define a~linear order $<_R$ on $\R\cup\R_{\mathrm{f}}$ and begin 
with nonnegative (false) decimals $d=\sum_{n\le k}d_n10^n$ and $e=\sum_{n\le k'}e_n10^n$, $d\ne e$. We set $d<_R e$ if and only if $d_i<e_i$ 
for the maximum $i\in\mathbb{Z}$ with $d_i\ne e_i$ (recall that $d_i=0$ for $i>k$ and $e_i=0$ for $i>k'$). If $d$ is negative and $e$ nonnegative, then $d<_R e$.
If $d$ and $e$ are distinct and both negative, we set $d<_R e$ if and only if $e_-<_R d_-$. We set $d\not<_R d$ 
for every $d\in\mathbb{R}\cup\mathbb{R}_{\mathrm{f}}$. Clearly, for $d,e\in\mathbb{R}_{\mathrm{t}}$ one has that $d<_R e$ if and only if $r(d)<r(e)$. Thus we have two linear orders,
$$
(\mathbb{R},\,<_R)\;\text{ and }\;(\mathbb{R}\cup\mathbb{R}_{\mathrm{f}},\,<_R)\;.
$$ 
Note that in the latter order, for any nonnegative $d\in\mathbb{R}_{\mathrm{t}}$ the pair $\overline{d}<_R d$ is a~{\em gap}, there is no element of $\mathbb{R}\cup\mathbb{R}_{\mathrm{f}}$ in-between, 
and the same holds for any negative $d\in\mathbb{R}_{\mathrm{t}}$ for the pair $d<_R\overline{d}$. One can think of the entries in these pairs as infinitesimally 
close, yet distinct.

Let $(X,<_X)$ be a~linear order. A~set $Y\subset X$, or a~sequence $(x_n)\subset X$, is {\em bounded from above} if there is an $x\in X$ such that $y\le_X x$
for every $y\in Y$, or $x_n\le_X x$ for every $n\in\mathbb{N}$. Any such $x$ is an {\em upper bound} of $Y$, or of $(x_n)$. Bounding from below and lower bounds
are defined in an analogous way. A~set or a~sequence is {\em bounded}, if it is bounded both from above and from below. 
We say that $x\in X$ is a~{\em supremum of $Y\subset X$}, and write $x=\sup(Y)$, if $x$ is the least upper bound of $Y$: $x$ is an
upper bound of $Y$ and if $x'\in X$ is any upper bound of $Y$ then $x\le_X x'$. Similarly, $x\in X$ is an {\em infimum of $Y\subset X$}, in symbols $x=\inf(Y)$, 
if $x$ is the largest lower bound of $Y$. Infimum and supremum, if they exist, are unique. We prove the key result that suprema and infima exist for sets of decimals, and also when false ones are added.

\begin{propo}\label{supremum}
Consider the linear orders $(\mathbb{R}\cup\mathbb{R}_{\mathrm{f}},<_R)$ and $(\mathbb{R},<_R)$. 
\begin{enumerate}
    \item In $(\mathbb{R}\cup\mathbb{R}_{\mathrm{f}},<_R)$, every nonempty set $U\subset\mathbb{R}\cup\mathbb{R}_{\mathrm{f}}$ that is bounded from above has supremum, and every nonempty 
    set $U\subset\mathbb{R}\cup\mathbb{R}_{\mathrm{f}}$ that is bounded from below has infimum.
    \item In $(\mathbb{R},<_R)$, every nonempty set $U\subset\mathbb{R}$ that is bounded from above has supremum, and every nonempty set $U\subset\mathbb{R}$ 
    that is bounded from below has infimum.
\end{enumerate}
\end{propo}
\proof
It is clear that $d=\inf(U)$ in $(\mathbb{R}\cup\mathbb{R}_{\mathrm{f}},<_R)$ if and only if $d_-=\sup(U_-)$ for the set
$$
U_-:=\{e_-\;|\;e\in U\}\;,
$$
and that the same holds in $(\mathbb{R},<_R)$.
In parts~1 and~2 it therefore suffices to prove only the first claim. Moreover, 
the first claim in part~1 implies the first claim in part~2. For suppose that $U\subset\mathbb{R}$ is nonempty and bounded from above in $(\mathbb{R},<_R)$. Thus $U$ is bounded 
from above in $(\mathbb{R}\cup\mathbb{R}_{\mathrm{f}},<_R)$ and we can take $d:=\sup(U)$ in $(\mathbb{R}\cup\mathbb{R}_{\mathrm{f}},<_R)$. If $d\in\mathbb{R}$ then
$d=\sup(U)$ also in $(\mathbb{R},<_R)$. If $d\in\mathbb{R}_{\mathrm{f}}$ then it is easy to see that $\overline{d}=\sup(U)$ in $(\mathbb{R},<_R)$. Thus it suffices 
to prove just the first claim in part~1.

Let $U\subset\mathbb{R}\cup\mathbb{R}_{\mathrm{f}}$ be nonempty and bounded from above in $(\mathbb{R}\cup\mathbb{R}_{\mathrm{f}},<_R)$. We may assume that
either every $d'\in U$ is positive, or that every $d'\in U$ is negative. We treat the former case first. Boundedness from above implies that orders of (false) decimals in $U$ 
attain the maximum value $k\in\mathbb{N}_0$. Let $U_k\subset U$ be the (false) decimals in $U$ with order $k$, and $d_k\in D$ be the maximum digit at $10^k$ in $d'\in U_k$. 
If the nonempty nested sets $U\supset U_k\supset U_{k-1}\supset\dots\supset U_l$ and digits $d_k,d_{k-1},\dots,d_l\in D$ are already defined, we set $U_{l-1}$ to be
the (false) decimals in $U_l$ with digit  $d_l$ at $10^l$, and $d_{l-1}\in D$ to be the maximum digit at $10^{l-1}$ in $d'\in U_{l-1}$. We set 
$$
d:=\sum_{n\le k}d_n10^n\in\mathbb{R}\cup\mathbb{R}_{\mathrm{f}}
$$
(clearly, $d>_R 0_R$) and show that $d=\sup(U)$.

Let $e=\sum_{n\le k'}e_n10^n\in U$ be arbitrary, so that $k'\le k$. If $k'<k$ then $e<_R d$. If $k'=k$ then $e_k\le d_k$. If $e_k<d_k$ then $e<_R d$. Else $k'=k$,
$e_k=d_k$, and we look for the minimum $l\in\mathbb{Z}$ such that $e_k=d_k$, $e_{k-1}=d_{k-1}$, $\dots$, $e_l=d_l$. If it does not exist then $e=d$. Else 
$e_{l-1}<d_{l-1}$ and $e<_R d$. Thus $d$ is an upper bound of $U$. Let $e=\sum_{n\le k'}e_n10^n\in\mathbb{R}\cup\mathbb{R}_{\mathrm{f}}$ be any (false) decimal with $e<_R d$, 
so that $k'\le k$. If $k'<k$, we take any $d'\in U_k$ and see that $e<_R d'$. If $k'=k$ then $e_k\le d_k$. If $e_k<d_k$, we take any $d'\in U_{k-1}$ and see that 
$e<_R d'$. If $k'=k$ and $e_k=d_k$ then we look for the minimum $l\in\mathbb{Z}$ such that $e_k=d_k$, $e_{k-1}=d_{k-1}$, $\dots$, $e_l=d_l$. Since $e<_R d$, the number 
$l$ exists and $e_{l-1}<d_{l-1}$. We take any $d'\in U_{l-2}$ and see that $e<_R d'$. Thus $d$ is the least upper bound of $U$ and $d=\sup(U)$.

We consider the latter case when every $d'\in U$ is negative. We define $k\in\mathbb{N}_0$ as the minimum order of (false) decimals in $U$. We define the nonempty nested 
sets $U\supset U_k\supset U_{k-1}\supset\dots$ and digits $d_k,d_{k-1},\dots\in D$ as in the former case, with the only change that the ``maximum digit at $10^{l-1}$''
is replaced with the ``minimum digit at $10^{l-1}$'', and set like before
$$
d:=-\sum_{n\le k}d_n10^n\in\mathbb{R}\cup\mathbb{R}_{\mathrm{t}}\;,
$$
thus $d<_R 0_R$. The proof that $d=\sup(U)$ goes as in the former case, with the only changes that in the previous paragraph we add two minus signs to the two (false) decimals $e$, and reverse each of the ten inequalities 
(listed in the order of reading) ``$k'\le k$'', ``$k'<k$'', ``$e_k\le d_k$'', ``$e_k<d_k$'', ``$e_{l-1}<d_{l-1}$'', ``$k'\le k$'', ``$k'<k$'', ``$e_k\le d_k$'', ``$e_k<d_k$'', 
and ``$e_{l-1}<d_{l-1}$''.
\eproof

A~{\em subsequence} $(x_{m_n})$ of a~sequence $(x_n)$ is determined by a~sequence of natural numbers $m_1<m_2<\dots$\;. If $(x_n)\subset X$
is a~sequence in a~linear order $(X,<_X)$, we say that it is {\em non-decreasing} if $x_1\le_X x_2\le_X\dots$\;. Non-increasing  sequences are defined 
in an analogous way. By a~{\em monotone sequence} we mean a~non-decreasing or a~non-increasing sequence. If $d\in\mathbb{R}$ is a~decimal and 
$(a_n)\subset\mathbb{Q}_D$ is a~sequence of fractions, we say that {\em $d$ is a~limit of $(a_n)$} if for every $k\in\mathbb{N}$ there is an $n_0\in\mathbb{N}$ such that for 
every integer $n\ge n_0$,
$$
|r(d\,|\,n)-a_n|<1/k\;.
$$
We prove that every appropriately bounded monotone sequence in $\mathbb{Q}_D$ has a~limit. 

\begin{propo}\label{zeroth_res}
Let $(X,<_X)$ be a~linear order.
\begin{enumerate}
    \item Every sequence $(x_n)\subset X$ has a~monotone subsequence.
    \item Every sequence $(a_n)\subset\mathbb{Q}_D$ that is non-decreasing and bounded from above has a~limit $d\in\mathbb{R}$. The same holds if $(a_n)$ is
     non-increasing and bounded  from below.
    \item For any two decimals $d<_R e$ there exist $k$ and $n_0$ such that for every integer $n\ge n_0$ one has that
    $r(e\,|\,n)-r(d\,|\,n)>1/k$.
\end{enumerate}
\end{propo}
\proof
1. Let $(x_n)$ be any sequence in $X$. We call $m\in\mathbb{N}$ a~{\em horizon} if $x_n<_X x_m$ for every $n>m$. If $(x_n)$ has infinitely many 
horizons, $m_1<m_2<\dots$, then $x_{m_1}>_X x_{m_2}>_X\dots$ is a~monotone subsequence. Else we take an $m_1\in\mathbb{N}$ larger than all
horizons in $(x_n)$. There is an $m_2>m_1$ such that $x_{m_1}\le_X x_{m_2}$. But there is also an $m_3>m_2$ such that $x_{m_2}\le_X x_{m_3}$, and so on. This
way we again obtain a~monotone subsequence.

2. By changing signs we may assume that we are given fractions $b$ and $a_n$ in $\mathbb{Q}_D$, $n\in\mathbb{N}$, such that 
$$
0\le a_1\le a_2\le\dots\le b\;\text{ or }\;a_1\ge a_2\ge\dots\ge 0\;.
$$
We treat the former case first. Using the first claim in part~2 of Proposition~\ref{supremum}, we set $d:=\sup_{<_R}(\{r^{-1}(a_n)\;|\;n\in\mathbb{N}\})$ 
and show that $d$ is a~limit of $(a_n)$. Let a~$k$ be given. We take an $m\in\mathbb{N}$ such that $2/10^m<1/k$ and set 
$$
e:=d\,|\,m-_R 1\cdot10^{-m}\in\mathbb{R}_{\mathrm{t}}\;. 
$$
It follows that for every integer $n\ge m$ we have the bound
\begin{equation}\label{the_bound}
r(d\,|\,n)-r(e)<2/10^m\;.    
\end{equation}

By the definition of supremum and by monotonicity of $(a_n)$ it follows that for every large enough $n\ge m$ we have that 
$$
e<_R r^{-1}(a_n)\le_R d<_R d\,|\,n+_R 1\cdot10^{-n}\;.
$$
Hence there is an $n_0$ with $n_0\ge m$ and such that for every integer $n\ge n_0$,
$$
e-_Rd\,|\,n<_R r^{-1}(a_n)-_Rd\,|\,n<_R 1\cdot10^{-n}\;.
$$
(We use the operation $-_R$ and the fact that $\mathbb{Q}$ is an ordered field.)
By applying the map $r$ and inequality (\ref{the_bound}) we get that for every integer $n\ge n_0$,
$$
-1/k<-2/10^m<a_n-r(d\,|\,n)<1/10^n<1/k\;.
$$
Thus $d$ is a~limit of $(a_n)$.

We treat the latter case when $a_1\ge a_2\ge\dots\ge 0$. Using the second claim in part~2 of Proposition~\ref{supremum} we set $d:=\inf_{<_R}(\{r^{-1}(a_n)\;|\;n\in\mathbb{N}\})$ and show that 
$d$ is a~limit of $(a_n)$. Let a~$k$ be given. We take an $m\in\mathbb{N}$ such that $1/10^m<1/k$ and set 
$$
e:=d\,|\,m+_R 1\cdot10^{-m}\in\mathbb{R}_{\mathrm{t}}\;. 
$$
It follows that for every integer $n\ge m$ we have the bound
\begin{equation}\label{the_bound2}
r(e)-r(d\,|\,n)<1/10^m\;.    
\end{equation}

By the definition of infimum and by monotonicity of $(a_n)$ it follows that for every large enough $n\ge m$ we have that 
$$
e>_R r^{-1}(a_n)\ge_R d\ge_R d\,|\,n\;.
$$
Like before we subtract $d\,|\,n$, apply the map $r$, use inequality (\ref{the_bound2}) and get that for every integer $n\ge n_0\ge m$,
$$
1/k>1/10^m>a_n-r(d\,|\,n)\ge0\;.
$$
Thus $d$ is a~limit of $(a_n)$.

3. We assume that $d$ and $e$ are nonnegative, other cases are reduced to this situation by the map $d\mapsto d_-$
and/or by considering the pairs of decimals $d<_R 0_R$ and $0_R<_R e$. We take maximum $n\in\mathbb{Z}$ such that $d_n\ne e_n$, then $d_n<e_n$. If $e_n-d_n\ge2$ then
$r(e\,|\,j)-r(d\,|\,j)>10^n$ for every $j\ge\max(1,-n)$. Let $e_n-d_n=1$ and $m\in\mathbb{Z}$, $m<n$, be maximum with $d_m\ne9$  (such $m$ exists by 
property (ii) of decimals). Then $r(e\,|\,j)-r(d\,|\,j)>10^m$ for every $j\ge\max(1,-m)$.
\eproof

After these preparations we introduce generating sequences of decimals.
We say that a~sequence $(a_n)=(a_1,a_2,\dots)\subset\mathbb{Q}_D$ is {\em Cauchy}
if for every $k$ there is an $n_0$ such that for all integers $m,n\ge n_0$,
$$
|a_m-a_n|<1/k\;.
$$
For example, for every $d\in\mathbb{R}$ the sequence of truncations $(r(d\,|\,n))$ is Cauchy. It is clear that every Cauchy sequence $(a_n)\subset\mathbb{Q}_D$
is bounded. Two Cauchy sequences $(a_n),(b_n)\subset\mathbb{Q}_D$ 
are {\em mutually close} if for every $k$ there is an $n_0$ such that for all integers $m,n\ge n_0$,
$$
|a_m-b_n|<1/k\;.
$$
Equivalently we may take $m=n$ here. 

The key definition in our approach is that a~decimal $d\in\mathbb{R}$ is {\em generated by a~Cauchy sequence $(a_n)\subset\mathbb{Q}_D$}, we also say 
that $(a_n)$ is a~{\em generating sequence of $d$}, if $d$ is a~limit of $(a_n)$: for every $k$ there is an $n_0$ such that
$$
n\ge n_0\Rightarrow |r(d\,|\,n)-a_n|<1/k\;.
$$
We write $d=G(a_n)$, as in \cite[Section~2]{klaz} (there $X=G(a_n)$ with $X\subset\mathbb{Q}$ being a~cut). We prove three simple but important 
results on the relation $d=G(a_n)$.

\begin{propo}\label{first_res}
For every decimal $d\in\mathbb{R}$ there exists a~Cauchy sequence $(a_n)\subset\mathbb{Q}_D$ such that $d=G(a_n)$.
For every Cauchy sequence $(a_n)\subset\mathbb{Q}_D$ there exists a~unique decimal $d\in\mathbb{R}$ such that $d=G(a_n)$.
\end{propo}
\proof
The first claim is immediate because for any $d\in\mathbb{R}$ one has that $$d=G(r(d\,|\,n))\;.$$ To prove the second claim, we assume that
$(a_n)\subset\mathbb{Q}_D$ is a~Cauchy sequence. By part~1 of Proposition~\ref{zeroth_res}, $(a_n)$ has a~bounded monotone subsequence $(a_{m_n})$. 
By part~2 of Proposition~\ref{zeroth_res}, $(a_{m_n})$ has a~limit $d\in\mathbb{R}$. Since the sequence $(a_n)$ is Cauchy, it follows from the inequality
$$
|r(d\,|\,n)-a_n|\le|r(d\,|\,n)-a_{m_n}|+|a_{m_n}-a_n|
$$
that $d=G(a_n)$. Uniqueness of $d$ follows from part~3 of Proposition~\ref{zeroth_res}.
\eproof

\begin{propo}\label{second_res}
If $d$ is a~decimal, $d=G(a_n)$ and also $d=G(b_n)$, then the Cauchy sequences $(a_n)$ and $(b_n)$ are mutually close.
\end{propo}
\proof
This follows from the inequality
$$|a_n-b_n|\le|a_n-r(d\,|\,n)|+|r(d\,|\,n)-b_n|\;.$$
\eproof

\begin{propo}\label{third_res}
If $d$ and $e$ are decimals, $d=G(a_n)$, $e=G(b_n)$ and the Cauchy sequences $(a_n)$ and $(b_n)$ are mutually close, then $d=e$.
\end{propo}
\proof
If $d$, $e$, $(a_n)$ and $(b_n)$ are as stated, then we easily see from an inequality similar to the previous one that the sequences $(r(d\,|\,n))$ 
and $(r(e\,|\,n))$ are mutually close. By part~3 of Proposition~\ref{zeroth_res}, $d=e$.
\eproof

\section{Dedekind's theorem}

We construct by means of decimals the {\em complete ordered field $\mathbb{R}$}. 
As in \cite{klaz} we recall that this is an algebraic structure
$$
F=(F,\,0_F,\,1_F,\,+_F,\,\cdot_F,\,<_F)
$$
such that $F$ is a~set, $0_F$ and $1_F$ are two different elements of $F$, 
$$
+_F,\,\cdot_F\colon F\times F\to F
$$ 
are two binary operations on $F$,
$<_F\,\subset F\times F$ is a~linear order on $F$, and the following axioms hold. Both operations $+_F$ and $\cdot_F$ are associative and 
commutative, and are bound by the distributive law. The element $0_F$ is neutral with respect to $+_R$, and $1_F$ with respect to $\cdot_R$. Every element
$x\in F$ has an (additive) inverse $-x$ to $+_R$, and every element $x\in F\setminus\{0_F\}$ has a~(multiplicative) inverse $x^{-1}$ to $\cdot_R$. 
The linear order $<_R$ relates to the operations by two {\em order-operation axioms}\,---\,for any $a,b,c\in F$,
$$
a<_F b\Rightarrow a+_F c<_F b+_F c\;\text{ and }\;0_F<_F a\wedge 0_F<_F b\Rightarrow 0_F<_F a\cdot_F b\;.
$$
These are the axioms of an {\em ordered field $F$}. Its {\em completeness} means that any nonempty set $A\subset F$ that is in $(F,<_F)$ bounded from above has supremum in $(F,<_F)$. 

We employ decimals and prove that a~complete ordered field exists. Compared to the proof by
Dedekind's cuts there is one more technical complication, in the ring $\mathbb{Q}_D$ one cannot divide. We fix it by approximate division.

\begin{lemma}\label{division}
For every nonzero fraction $a\in\mathbb{Q}$ and every number $k$ there exists a~fraction $b\in\mathbb{Q}_D$ such that
$$
|ab-1|<1/k\;.
$$
\end{lemma}
\proof
Let a~nonzero $a\in\mathbb{Q}$ and a~number $k$ be given.
We write $a=\frac{m}{n}$ for a~nonzero integer $m$ and an $n\in\mathbb{N}$, and take an $l\in\mathbb{N}$ such that $|m|/10^l<1/k$. We take an integer
$p$ such that $|pm-10^l|<|m|$. Then the fraction $b:=pn/10^l\in\mathbb{Q}_D$ satisfies that
$$
|ab-1|=\bigg|\frac{mpn}{n10^l}-1\bigg|=\bigg|\frac{mp-10^l}{10^l}\bigg|<\frac{|m|}{10^l}<\frac{1}{k}\;.
$$
\eproof

Now we can prove by decimals that there exists a~complete ordered field. This was first proven  by means of cuts by R.~Dedekind, and published by him in 
\cite{dede}. We hasten to add that our formulation of his result is anachronistic. 

\begin{thm}[R.~Dedekind, 1858]\label{dede}
There exists a~complete ordered field $\R$.
\end{thm}
\proof (Proof by decimals.)
We start from the ordered field of fractions $$\mathbb{Q}=(\mathbb{Q},\,0,\,1,\,+,\,\cdot,\,<)$$ which we take for granted\,---\,we have been working with it already\,---\,and construct from it a~complete ordered 
field $$\R=(\mathbb{R},\,0_R,\,1_R,\,+_R,\,\cdot_R,\,<_R)\;.$$ We use representations of decimals as $d=G(a_n)$, guaranteed by Proposition~\ref{first_res}. For $a\in\mathbb{Q}_D$ let $c_a=(a,a,\dots)\subset\mathbb{Q}_D$ be the constant (Cauchy) sequence of $a$s. The set $\mathbb{R}$ was already 
defined as the set of all decimals, also $0_R=G(c_0)$ was already defined, and we set $1_R:=1\cdot10^0\in\mathbb{R}_{\mathrm{t}}$, so that $1_R=G(c_1)$. We define for any two decimals $d=G(a_n)\in\mathbb{R}$ and 
$e=G(b_n)\in\mathbb{R}$ that
$$
d+_R e:=G(a_n+b_n),\ d\cdot_R e:=G(a_n\cdot b_n)\;\text{ and }\;-d:=G(-a_n)\;.
$$
These definitions are correct because for any two Cauchy sequences $(a_n),(b_n)\subset\mathbb{Q}_D$, the sequences $(a_n+b_n)$, $(a_n\cdot b_n)$ and $(-a_n)$
are contained in $\mathbb{Q}_D$ and are Cauchy. The definitions are also correct because the defined elements of $\mathbb{\R}$ do not depend
on the choice of generating sequences of decimals $d$ and $e$. We show it only for the addition $+_R$. If also $d=G(a_n')$ and $e=G(b_n')$ then by 
Proposition~\ref{second_res} the sequences $(a_n)$ and $(a_n')$, and $(b_n)$ and $(b_n')$, are mutually close. Thus also $(a_n+b_n)$ and $(a_n'+b_n')$  
are mutually close and by Proposition~\ref{third_res}, $G(a_n+b_n)=G(a_n'+b_n')$ as decimals. If we represent any $d\in\mathbb{R}_{\mathrm{t}}$ as $d=G(c_{r(d)})$, we see that 
the restrictions of the operations $+_R$ and $-_R$ (where $d-_R e:=d+_R (-e)$) to $\mathbb{R}_{\mathrm{t}}$ coincide with the earlier operations $+_R$ and $-_R$.

We define the multiplicative inverse to $d=G(a_n)$ if $d\ne0_R$. We may assume that $a_n\ne0$ for every $n$ and using Lemma~\ref{division} we take for every 
$n\in\mathbb{N}$ a~$b_n'\in\mathbb{Q}_D$ such that $|a_nb_n'-1|<1/n$. It is not hard to show that $(b_n')\subset\mathbb{Q}_D$ is a~Cauchy sequence and therefore (by Proposition~\ref{first_res}) we may set
$$
d^{-1}:=G(b_n')\;.
$$
It is also not hard to show that, like earlier, $d^{-1}\in\mathbb{R}$ does not depend on the choice of generating sequences $(a_n)$ and $(b_n')$. The ordered field $\mathbb{R}$
has been defined, since we defined $<_R$ earlier, and we check that the fifteen or so axioms of a~complete ordered field are satisfied in $\mathbb{R}$.

Associativity and commutativity of $+_R$ and $\cdot_R$ and the distributive law transfer straightforwardly from the field $\mathbb{Q}$, more precisely from the ring $\mathbb{Q}_D$. For example, $+_R$
is associative because addition is associative in $\mathbb{Q}_D$:
$$
(d+_R e)+_R f=G((a_n+b_n)+c_n)=G(a_n+(b_n+c_n)=d+_R(e+_R f)\;.
$$
In a~similar way we easily check the axioms for neutral and inverse elements: $d+_R 0_R=G(a_n+0)=G(a_n)=d$, 
$d\cdot_R 1_R=G(a_n\cdot1)=G(a_n)=d$, $d+_R(-d)=G(a_n+(-a_n))=G(c_0)=0_R$, and for $d\ne0_R$ also $d\cdot_R d^{-1}=G(a_n\cdot b_n')=G(c_1)=1_R$ by 
Proposition~\ref{third_res} because the Cauchy sequences $(a_n\cdot b_n')$ and $c_1$ are mutually close. 

We already know that $<_R$ is a~linear order on $\mathbb{R}$ and that every nonempty subset of $\mathbb{R}$ that is bounded from above has supremum in $(\mathbb{R},<_R)$, this is the first claim
in part~2 of Proposition~\ref{supremum}. It only remains to check the order-operation axioms. It is easy to see that if $d=G(a_n)$ and $e=G(b_n)$ 
are decimals then $d<_R e$ if and only if there exist numbers $k,n_0$ such that for every integer $n\ge n_0$ one has that $b_n-a_n>1/k$. Let $f=G(c_n)$
be another decimal and let $d<_R e$, witnessed by some numbers $k,n_0$. Then
$$
d+_R f=G(a_n+c_n)<_R e+_R f=G(b_n+c_n)
$$
because $(b_n+c_n)-(a_n+c_n)=b_n-a_n>1/k$ for every $n\ge n_0$. If $0_R<_R d,e$ and numbers $k$ and $n_0$ are such that $a_n,b_n>1/k$  for every $n\ge n_0$, then
$$
0_R <_R d\cdot_R e=G(a_n\cdot b_n)
$$
because $a_nb_n>1/k^2$  for every $n\ge n_0$. The proof that $\mathbb{R}$ is a~complete ordered field is complete.
\eproof

\noindent
One can also prove that the ordered field $\mathbb{Q}$ embeds densely in $\mathbb{R}$, and that every complete ordered field is 
isomorphic to $\mathbb{R}$ as an ordered field, but we will not go in these results here.

A~possible objection may be that the previous proof just restates in a cumbersome way the well-known Cantor's construction
of the real numbers in terms of equivalence classes of Cauchy sequences of fractions. Our answer is that the substantial difference\,---\,as far as we know, ignored up to now\,---\,between the decimal (or Dedekind's) construction on one side and Cantor's constructions on the other side is that any decimal $d$ (or any cut)
is a~countable set, or in fact even a~hereditarily at most countable set (which means that additionally every element of $d$ is at most countable, 
every element of every element of $d$ is at most countable, and so on), whereas any equivalence class in Cantor's construction is an uncountable set. See \cite[Section~4]{klaz} for more discussion, here our interest lies in computing the arithmetic operations on decimals. 

\section{Arithmetic of decimals}

The arithmetic operations $+_R$ and $\cdot_R$ on decimals, defined in the previous proof, 
are troublesome from the perspective of digits. Consider the next computation (from now on we use for brevity standard notation for decimals):
$$
d+_R e:=0.33333\dots+_R 0.66666\dots=\left\{
\begin{array}{ll}
  f:=1.00000\dots   & \text{or} \\
  g:=0.99999\dots   & ?
\end{array}\right.
$$
If the ellipses in $d$ and $e$ mean that after the decimal point there is only the digit $3$, resp. $6$, then the 
correct sum is the terminating decimal $f=1=1\cdot10^0$, and not the false decimal $g=0.99999\dots\;$. However, sums of truncations of $d$ and $e$ 
formally converge to $g$.

But the ellipses in $d$ and $e$ might also mean that in $d$ after a~long segment of $3$s other digits appear, and similarly in $e$. 
Then we may get the sum $d+_R e$ in both forms, it may be either $f$ with some nonzero digits after a~long segment of $0$s, or $g$ with some
digits less than nine after a~long segment of $9$s. For example,
$$
0.3333371\dots+_R 0.6666640\dots=1.0000011(2)\dots
$$
or
$$
\;\;0.3333324\dots+_R 0.6666654\dots=0.9999978(9)\dots\;,
$$
always with two possibilities in the last digit. We assume that the usual arithmetic of decimals coincides with the 
arithmetic defined in the proof of Theorem~\ref{dede}, which is true but should be proven, and return to this in Theorem~\ref{thm_weak_comp_addi}.
The fact that the leading digit at $10^0$ in $d+_R e$ (resp. $d\cdot_R e$) is sensitive to arbitrarily distant changes of digits in $d$ and $e$ causes the
{\em incomputability} of the operation $+_R$ (resp. $\cdot_R$). Later we show that, on the other hand, $+_R$ 
and $\cdot_R$ can be {\em weakly computed}.

Before we formalize these kinds of computation, we describe for contrast the arithmetic of the classical number system $\mathbb{Q}_p$ of {\em $p$-adic numbers}, associated to  a~prime number $p\in\mathbb{N}$. Its addition and  multiplication {\em are} computable. For more information on $\mathbb{Q}_p$ see, for example, \cite{robe}. 
Elements of $\mathbb{Q}_p$ are quite similar to decimals, a~$p$-adic number~$a$
is a~formal infinite linear combination
$$
a=\sum_{n\ge k}a_np^n
$$
such that $k,n\in\mathbb{Z}$ with $k\le 0$, and that  for every $n$ one has $a_n\in\{0,1,\dots,p-1\}=:P_p$. Like for decimals, we call $k$ the {\em order of $a$}  
and require that $a_k=0\Rightarrow k=0$. We also call each $a_n$ the {\em digit at $p^n$ (in $a$)} and set the $a_n$ with $n<k$ to zero. Two $p$-adic numbers are
equal iff they have equal orders and equal digits at $p^n$ for every $n\in\mathbb{Z}$. 

These numbers are simpler than decimals, there are no signs and the only restriction on digits is the analogue of (i). 
The addition $+_p$ and the multiplication $\cdot_p$ of $p$-adic numbers are also simpler than for decimals. We describe in detail for $a,b\in\mathbb{Q}_p$ their sum
$$
a+_p b=\sum_{n\ge k}a_np^n+_p \sum_{n\ge k'}b_np^n=\sum_{n\ge k''}c_np^n\in\mathbb{Q}_p\;.
$$
First we set $k_0:=\min(k,k')$, let $c_{k_0}\in P_p$ be the unique number congruent to $a_{k_0}+b_{k_0}$ modulo $p$, and let $t_{k_0}\in\{0,1\}$ be $1$ 
iff $a_{k_0}+b_{k_0}\ge p$ (and $0$ else). If the digits  $c_{k_0},c_{k_0+1},\dots,c_l\in P_p$ and carries $t_{k_0},t_{k_0+1},\dots,t_l\in\{0,1\}$
are already defined, we set $c_{l+1}\in P_p$ to be the unique number congruent to $a_{l+1}+b_{l+1}+t_l$ modulo $p$, and set $t_{l+1}\in\{0,1\}$ to be $1$ 
iff $a_{l+1}+b_{l+1}+t_l\ge p$ (and $0$ else). The order $k''$ of the sum is then the minimum $j\le0$ with $c_j\ne0$ if such $j$ exists, and $k''=0$ else. This is the well-known term-wise addition modulo $p$ with carries to higher orders. 

The multiplication $\cdot_p$ is defined similarly, although in a~little more complicated way, via distributivity 
and again with carries to higher orders. It is well known that 
$(\mathbb{Q}_p,0_p,1_p,+_p,\cdot_p)$, where $0_p:=\sum_{n\ge0}0p^n$ and $1_p:=1p^0+_p \sum_{n\ge1}0p^n$, is a~field. Its main property is completeness with respect to a~certain norm, but this is not important here. Here the important aspect of $\mathbb{Q}_p$ is that, in contrast 
with $+_R$ and $\cdot_R$, both operations $+_p$ and $\cdot_p$ are computable. We introduce this notion in a~general setting.

Let $A\ne\emptyset$ be  a~finite alphabet and let
$$
A^{\infty}:=A^{\mathbb{N}_0}=\{u=u_0u_1u_2\dots\;|\;u_n\in A\}
$$
be the set of one-way infinite words over $A$. 
An~{\em encoding} is any nonempty subset $X\subset A^{\infty}$. An~{\em operation on $X$} is any map
$$
F\colon X\times X\to X\;.
$$
For example, $F$ may be addition or multiplication of $p$-adic numbers or of decimals, given in some encoding. We say that $F$ is {\em computable} if there exists a~machine $M$ computing~$F$, and that $F$ is {\em weakly computable} if there exists
a~machine $M$ weakly computing~$F$. We introduce machines and computing; weak computing will be defined later.

Suppose that $X\subset A^{\infty}$ is an encoding.
A~{\em machine $M$} is any Turing machine, with a~finite {\em set of states} that contains an {\em initial state} and a~{\em terminal state}, and with a~finite {\em alphabet $A'$} such that $A'\supset A\cup\{0,1,\epsilon\}$, where $\epsilon\in A'\setminus A$ is the {\em empty cell symbol}. $M$ operates on three biinfinite tapes $T_1,T_2,T_3$\,---\,each consists of cells $c$
labeled by the integers $\mathbb{Z}$\,---\,via three corresponding reading/writing heads $H_1,H_2,H_3$. Tapes $T_1$ and $T_2$ are input tapes, and $T_3$ is an input-output tape. 

$M$ works as follows. Its inputs are the triples 
$$
(x,\,y,\,m)\in X\times X\times\N_0\;,
$$ 
where the word $x=x_0x_1\dots$ is stored on $T_1$ so that the letter $x_n$ is stored in cell $n$ and each cell $n<0$ stores $\epsilon$, the word $y$ is stored on $T_2$
in the same way, and the number $m$ is stored on $T_3$ in the binary form, with the least significant digit at the left in cell $0$ (we store binary numbers on tapes in the order opposite to the reading order). For example, if $m=16$ then
$$
T_3=\ds\,\epsilon\,\epsilon\,\underline{0}\,0\,0\,0\,1\,\epsilon\,\epsilon\,\ds\;,
$$
with the (content of the) cell $0$ underlined; we will mark the cell $0$ in this way. Heads $H_1$ and $H_2$ can only read, but $H_3$ can read and write. The {\em computation (run) of $M$ on the input $(x,y,m)$} starts with $M$ in the initial state and with each head 
$H_i$, $i=1,2,3$, positioned at the cell $c_i=0$ of $T_i$. The action of $M$ in one step is fully determined by its current state and by the content of the cells $c_1,c_2,c_3$. 
The machine $M$ consults its (finite) instruction table and does this: if $M$ is in the final state then the computation terminates, else $M$ rewrites $c_3$ to the prescribed new (possibly the same) symbol in $A'$, moves each head $H_1,H_2,H_3$ in the prescribed way by $\pm1$ or $0$ position, and jumps in the prescribed new (possibly the same) state. 

\begin{defi}\label{defi_comput}
An operation $F\colon X\times X\to X$
 on an encoding $X$ is {\em computed by a~machine $M$} if for every input $(x,y,m)\in X^2\times\N_0$ the computation of $M$ on it terminates and the content of the cell $0$ of $T_3$ is $z_m$ (other cells of $T_3$ are filled with $\epsilon$), where $F(x,y)=z_0z_1\ds\;$\,---$M$ computes the $m$-th letter of the result of the operation $F$ on the words $x,y$. We say that $F$ is {\em computable} if there exists a~machine $M$ computing $F$.
 \end{defi}

It is clear that addition and multiplication of $p$-adic numbers, of which we described the former in detail, are computable in this sense. We state it formally. We set $A:=P_p\cup\{\xi\}$, where $\xi$
is a~separation symbol, and encode any $p$-adic number $a=\sum_{n\ge k}a_np^n$ as a~word $u(a)$ in $A^{\infty}$ so that $u(a)$ starts with the binary digits of $|k|$, then the symbol $\xi$ follows, and $u(a)$ is completed by the $p$-adic digits $a_k,a_{k+1},\dots$ of $a$ in this order. This gives the {\em canonical encoding of $p$-adic numbers}
$$
X_p:=\{u(a)\;|\;a\in\mathbb{Q}_p\}\;.
$$
We regard $+_p$ and $\cdot_p$ as operations on $X_p$. 
For example, if $a\in\mathbb{Q}_p$ has order $0$ and $b\in\mathbb{Q}_p$ has order $-6$,
then the words $u(a)$ and $u(b)$ appear on a~tape $T$ as
$$
T=\dots\epsilon\,\epsilon\,\underline{0}\,\xi\,a_0\,a_1\,\dots\;
\text{ and }\;T=\dots\epsilon\,\epsilon\,\underline{0}\,1\,1\,\xi\,b_{-6}\,b_{-5}\,\dots\;,
$$
respectively.

\begin{thm}\label{thm_comp}
For any prime number $p$, both operations $+_p$ and $\cdot_p$ on the canonical encoding $X_p$ of $p$-adic numbers are computable. 
\end{thm}
\proof
The instruction table of a~machine $M$ computing $+_p$ is outlined on a~high level in the above description of the sum $a+_p b$. For 
the product $a\cdot_p b$ we proceed similarly. As common in design of algorithms, we trust in our intuition that one can convert such
high level descriptions in actual instruction tables of Turing machines.
\eproof

\noindent
The arithmetic of $p$-adic numbers is even nicer than that: the digit at $p^n$ in a~sum or in a~product depends only on the digits at $p^r$ of the 
operands such that $r\le n$. 

We show that the operations $+_R$ and $\cdot_R$ on decimals are not computable. For decimals we set $A:=D\cup\{-,\xi\}$, where $D=\{0,1,\ds,9\}$, $\xi$ has the previous meaning and $-$ is the possible minus sign of the encoded decimal, and 
encode a~decimal $d\in\mathbb{R}$ as a~word $u(d)\in A^{\infty}$ like in the case of $p$-adic numbers, with the only difference that $u(d)$ may start with the minus sign and that the indices 
of digits decrease. This is the {\em canonical encoding of decimals} 
$$
X_R:=\{u(d)\;|\;d\in\mathbb{R}\}\;.
$$
We again regard $+_R$ and $\cdot_R$ as operations on $X_R$. For example, if $d\in\mathbb{R}$ is negative and has order $2$, then the word $u(d)$ appears on a~tape $T$ as
$$
T=\dots\,\epsilon\,\epsilon\,\underline{-}\,0\,1\,\xi\, d_2\, d_1\dots\;.
$$
In particular, any nonnegative decimal $d$ with order $0$ has its digit $d_0$ stored in cell $2$.

\begin{thm}\label{thm_incomp}
Neither of the operations $+_R$ and $\cdot_R$ on $X_R$ is computable.
\end{thm}
\proof
We prove it in detail for addition, for multiplication one argues similarly. Suppose in the way of contradiction that $M$ is a~machine computing the operation $+_R$ on $X_R$. We consider the decimals $x:=0.66666\dots$ and $y:=0.33333\dots$ of the earlier example. First we run $M$ on the input
$$
(u(x),\,u(y),\,2)=(0\,\xi\,0\,6\,6\,6\,6\,6\,\dots,\,0\,\xi\,0\,3\,3\,3\,3\,3\,\dots,\,2)\;.
$$
After some $r\in\N$ steps the computation of $M$ terminates. Then the tape $T_3$ is
$$
T_3=\dots\,\epsilon\,\epsilon\,\underline{1}\,\epsilon\,\epsilon\,\dots\;,
$$
because $u(x+_R y)=u(1)=0\xi 100\dots\;$. We change $x$ to $x'$ by changing in $x$ one digit, the digit $6$ at $10^{-r}$ is replaced by $5$. Then
$$
x'+_R y=0.\underbrace{99\dots 9}_{r}00\dots\;.
$$
We let $M$ run again, now on the input $(u(x'),u(y),2)$. During the first $r$ steps of the computation, the machine $M$ cannot detect that the input has changed 
because the head $H_1$ starts moving from the cell $0$ and therefore $M$ works completely identically to the previous run. It terminates again after $r$ steps and also the output content of $T_3$ is the same,
$T_3=\dots\epsilon\epsilon\underline{1}\epsilon\epsilon\dots\;$. But this contradicts the assumption that $M$ computes $+_R$ because now the output should be
$$
T_3=\dots\,\epsilon\,\epsilon\,\underline{0}\,\epsilon\,\epsilon\,\dots\;.
$$

The argument for incomputability of the product $\cdot_R$ is similar, one may for example use the equality
$$
0.33333\ldots\,\cdot_R 3=1\;.
$$
\eproof

The canonical encoding $X_R$ is just one of many possible encodings of decimals. The previous theorem should apply to all those of them that are
computably equivalent to $X_R$; we formalize this notion. Suppose that $A_1$ and $A_2$ are nonempty finite alphabets 
and $X_i\subset A_i^{\infty}$, $i=1,2$, are two encodings. We say that $X_1$ and $X_2$ are {\em equivalent} if there exists a~bijection $f\colon X_1\to X_2$ 
such that both $f$ and its inverse $f^{-1}$ are 
computed by a~unary machine; $f$ is then called an {\em equivalence of $X_1$ and $X_2$}. A~{\em unary machine $M$} that {\em computes a~map $f\colon X_1\to X_2$} has alphabet containing $\{\epsilon,0,1\}\cup A_1\cup A_2$ 
and works as an ordinary machine, with the only difference that it has unary input and thus operates on only two tapes, an input tape $T_1$ and an input-output tape $T_2$. 
The input $(x,m)\in X_1\times\N_0$ is stored on $T_1$ and $T_2$ as before, and the output $T_2=\ds\epsilon\epsilon\underline{o}\epsilon\epsilon\ds$ means that $o=y_m$, where $f(x)=y_0y_1\ds\in X_2$. In all other respects unary machines function
as ordinary machines. 

If $F\colon X_1\times X_1\to X_1$ is an operation on $X_1$ and $g\colon X_1\to X_2$ is a~bijection, then by the {\em $g$-conjugate of $F$} 
we mean the operation $F'\colon X_2\times X_2\to X_2$ on $X_2$, defined by
$$
F'(x,\,y):=g(F(g^{-1}(x),\,g^{-1}(y)))\;.
$$

\begin{corollary}\label{cor_gen_coding}
For any equivalence of encodings $g\colon X_R\to X$, the $g$-conjugates of the operations $+_R$ and $\cdot_R$ are not computable.
\end{corollary}
\proof (Sketch.)
Let $g\colon X_R\to X$ be as stated and let $N$ and $N'$ be unary machines computing
$g$ and $g^{-1}$, respectively. We assume for the contrary that $M$ is a~machine computing the $g$-conjugate of $+_R$ (or of $\cdot_R$).
The composite machine 
$$
N'\circ M\circ N=N'(M(N))
$$ 
then computes the operation $+_R$ (or $\cdot_R$) on $X_R$, which contradicts Theorem~\ref{thm_incomp}. We leave to the interested reader to work out formal details of composition of machines and unary machines.
\eproof

We mention another encoding of decimals in use that is equivalent to $X_R$, the {\em scientific encoding (notation)} 
$$
X_S=\{v(d)\;|\;d\in\mathbb{R}\}\subset(D\cup\{-,\,\xi\})^{\infty}\;. 
$$
In it we encode zero as $v(0_R):=\xi-00\ds\;$. For a~nonzero decimal $d=(-)\sum_{n\le k}d_n10^n\in\mathbb{R}$
let the integer $m$ be maximum with $d_m>0$. Then $v(d)$ begins with the possible minus sign of $d$, followed by $\xi$, followed by the possible minus sign of $m$,
followed by the binary digits of $|m|$, followed by $\xi$, and $v(d)$ concludes with the digits $d_m,d_{m-1},\ds$ in this order.
For example, the estimate $1.7566\times10^{-15}$ m for the diameter of nucleus of the hydrogen atom, mentioned in Wikipedia, is in the encoding $X_S$ stored on a~tape $T$ as
$$
T=\ds\,\epsilon\,\epsilon\,\underline{\xi}-1\,1\,1\,1\,\xi\,1\,7\,5\,6\,6\,0\,0\,\ds\;.
$$
We leave to the interested reader to work out an equivalence between $X_S$ and $X_R$. One has to be careful, though, not to encode
$0_R$ as, say, $v(0_R):=\xi 0\xi 00\ds\;$. Then it would be possible for a~unary machine $M$ to decide whether the input decimal $d$ equals $0_R$ or not; $M$ would just check the digit after the second $\xi$. But such an encoding would not be equivalent with $X_R$ where for any given decimal $e\in\R$ and the input decimal $d$ the question ``does $d=e$?'' is algorithmicly undecidable.

Theorem~\ref{thm_incomp} and Corollary~\ref{cor_gen_coding} may appear a~little paradoxical because in practice we easily compute the sum or product of any pair of decimals\,---\,for example, in the proof of Theorem~\ref{thm_incomp} we computed the 
sum $x'+_R y$\,---\,and algorithms for $+_R$ and $\cdot_R$ appear in the literature. In these algorithms one freely jumps in a~decimal over any distance, 
indeed over its whole expansion, which no Turing machine can do. We formalize this kind of computation as the weak computation.

Suppose that $X\subset A^{\infty}$ is an encoding and that $F\colon X\times X\to X$ is an operation on $X$. 

\begin{defi}\label{defi_weak_comp}
We say that a~machine $M$ {\em weakly computes $F$} if for every 
pair $(x,y)\in X^2$ there exists a~number $h=h(x,y)\in\mathbb{N}$ such that the {\em enhanced machine $M_h$} computes as before for every
$m\in\N_0$ the $m$-th letter $z_m$ of $F(x,y)=z_0z_1\ds\;$. The enhancement means that the input for $M_h$ is $x$ on $T_1$ and $y$ on $T_2$ as before, but now also $m$ {\em and $h$} (both in binary notation) on $T_3$. The number $h$ is an additional {\em finite} piece of information which
helps $M$ to compute $F(x,y)$. We say that $F$ is {\em weakly computable} if a~machine $M$ weakly computes $F$.
\end{defi}
Weak computation of a~map $f\colon X_1\to X_2$ between two encodings by a~unary machine $M$ is defined in an analogous way. 

In the next two theorems we describe in reasonable detail machines $M$ weakly computing operations $+_R$ and $\cdot_R$ on $X_R$. As far as we know,
these are the first formal descriptions of algorithms for operations $+_R$ and $\cdot_R$ on decimals.

\begin{thm}\label{thm_weak_comp_addi}
The operation $+_R$ on $X_R$ is weakly computable.
\end{thm}
\proof
Let $d=(-)\sum_{n\le k}d_n10^n$ and $e=(-)\sum_{n\le k'}e_n10^n$ be two decimals. By the proof of Theorem~\ref{dede}, their sum
$f:=d+_R e\in\mathbb{R}$ is the limit of the Cauchy sequence
$$
(r(d\,|\,n)+r(e\,|\,n))\subset\mathbb{Q}_D\;.
$$
We weakly compute the sign of $f$, the order $k''\in\mathbb{N}_0$ of $f$, and for any $n\in\mathbb{Z}$ the digit $f_n$ at $10^n$ in $f$.
We can reduce the situation of the pair of decimals $d,e$ to two cases: $d, e\ge_R 0_R$ and $d\ge_R e_-\ge_R 0_R$, then always $f\ge_R 0_R$. We explain this reduction and the whole 
weak computation of $f$ at the end of the proof. Now we look at the digits $f_n$ and first treat the former case.

Let $d, e\ge_R 0_R$. The order $k''$ of $f$ is always encoded in an $h\in\mathbb{N}$ for the enhanced machine $M_h$. If $f\in\mathbb{R}_{\mathrm{t}}$ then the whole $f$ is
encoded in $h$. Suppose that $f\not\in\mathbb{R}_{\mathrm{t}}$. To find the digit $f_n$ we look for the maximum $l\in\mathbb{N}_0$ such that $d_{n-i}+e_{n-i}=9$ for 
$i=1,2,\dots,l$ (if $d_{n-1}+e_{n-1}\ne9$, we set $l=0$). Since $f$ is non-terminating, such $l$ exists. Then $f_n=d_n+e_n$  modulo $10$ if $d_{n-l-1}+e_{n-l-1}<9$, and $f_n=d_n+e_n+1$ modulo $10$ if $d_{n-l-1}+e_{n-l-1}>9$. 

Let $d\ge_R e_-\ge_R 0_R$. The order $k''$ of $f$ is always encoded in an $h\in\mathbb{N}$ for the enhanced machine $M_h$. Again,  
if $f\in\mathbb{R}_{\mathrm{t}}$ then the whole $f$ is encoded in $h$, and we assume that $f\not\in\mathbb{R}_{\mathrm{t}}$. To find the digit $f_n$, we look for the maximum 
$l\in\mathbb{N}_0$ such that $d_{n-i}=e_{n-i}$ for $i=1,2,\dots,l$ (if $d_{n-1}\ne e_{n-1}$, we set $l=0$). Since $f$ is non-terminating, such $l$ exists. Then $f_n=d_n-e_n$ modulo $10$ if 
$d_{n-l-1}>e_{n-l-1}$, and $f_n=d_n-e_n-1$ modulo $10$ if $d_{n-l-1}<e_{n-l-1}$. 

Any required information that a~machine $M$ attempting to compute $d+_R e$ cannot obtain alone may be encoded for the enhanced machine $M_h$ in $h\in\mathbb{N}$, as long as this information is finite. Thus we set, for example,
$$
h:=(2k''-1)\cdot2^r\;,
$$
where $k''\in\mathbb{N}_0$ is the order of $f$ and $r\in\mathbb{N}_0$ encodes in some way the fact whether $f$ is terminating, and if it is then 
$r$ encodes the whole $f\in\mathbb{R}_{\mathrm{t}}$.

To summarize, for a~given input $(u(d),u(e),m)$ the enhanced machine $M_h$ computes the $m$-th letter $z_m$ of $u(f)=u(d)+_R\, u(e)$ 
as follows. If $z_m$ corresponds to the part of $u(f)$ encoding the order of $f$, then $M_h$ recovers $z_m$ from $h$. The same, if 
$f$ is terminating. Suppose that $f\not\in\mathbb{R}_{\mathrm{t}}$ (which $M_h$ learns from $h$) and that $z_m$ is the sign of $f$ or some digit $f_n$. Then $d\ne e_-$, 
$M$ alone determines the sign of $f$ and reduces $d,e$ to one of the above two cases (simply by changing the signs of $d$ and $e$ as needed), and $M$ alone easily computes the digit $f_n$
by the above described procedures.
\eproof

\begin{thm}\label{thm_weak_comp_mult}
The operation $\cdot_R$ on $X_R$ is weakly computable.
\end{thm}
\proof
Let $d=(-)\sum_{n\le k}d_n10^n$ and $e=(-)\sum_{n\le k'}e_n10^n$ be two decimals. Now $f:=d\cdot_R e\in\mathbb{R}$ is the limit of the Cauchy sequence
$$
(r(d\,|\,n)\cdot r(e\,|\,n))\subset\mathbb{Q}_D\;.
$$
We weakly compute the sign of $f$, the order $k''\in\mathbb{N}_0$ of $f$, and for any given $n\in\mathbb{Z}$ the digit $f_n$ at $10^n$ in $f$.
We may assume that $d,e\ge_R 0_R$, then always $f\ge_R 0_R$. We explain this reduction and the whole weak computation of $f$ at the end of the proof.

So let $d, e\ge_R 0_R$. The order $k''$ of $f$ is always encoded in $h\in\mathbb{N}$ for the enhanced machine $M_h$. 
To find the digits $f_n$,  we consider for $l=1,2,\dots$ the terminating decimals
$$
f(l)=\sum_{n=-2l}^{k''}f(l)_n10^n:=(d\,|\,l)\cdot_R(e\,|\,l)\;.
$$
It is easy to see that for every $n\in\mathbb{Z}$ one has that $f(1)_n\le f(2)_n\le\dots\;$. So this sequence eventually stabilizes on a~value $g_n\in D$: for some $l_n\in\mathbb{N}$,  
$f(l_n)_n=f(l_n+1)_n=\dots=g_n$. If $f$ is non-terminating then $g_n=f_n$, but for $f\in\mathbb{R}_{\mathrm{t}}$ we may have 
$g_n=(\overline{f})_n$ (see the end of the proof of Theorem~\ref{thm_incomp}). This is resolved by enhanced machine.
To get $f_n$ we need an effective upper bound on $l_n$. Let $K:=\max(k,k')$ ($k$ and $k'$ is the order of $d$ and $e$, respectively). We claim that
$$
l_n\le\max(1,\,K-n+2)=:L_n\;.
$$
Suppose for the contrary that $f(l+1)_n-f(l)_n\ge1$ for some $l\in\mathbb{N}$ with $l\ge K-n+2$. Then
\begin{equation}\label{contr}
    K-l-n+1\le-1\;\text{ and }\;f(l+1)-f(l)\ge_R 10^n\;.
\end{equation}

From $10^{K+1}>_R d\ge_R d\,|\,l>_R d-_R 10^{-l}$ and analogous inequalities for $e$ we get that
$$
f\ge_R f(l+1)>_Rf(l)>_R f-_R 2\cdot_R 10^{K+1}\cdot_R 10^{-l}\;.
$$
From this we get, also using the second inequality in (\ref{contr}), that
$$
10^n\le_R f(l+1)-f(l)<_R
2\cdot10^{K+1}\cdot_R 10^{-l}\;\text{ and hence }\;K-l-n+1\ge0\;,
$$
which contradicts the first inequality in (\ref{contr}).

For given $d,e\in\mathbb{R}$ we set
$$
h:=(2k''-1)\cdot2^r\;,
$$
where $k''\in\mathbb{N}_0$ is the order of $f$ and $r\in\mathbb{N}_0$ encodes in some way the fact whether $f$ is a~terminating decimal, and 
if it is then $r$ encodes the whole $f\in\mathbb{R}_{\mathrm{t}}$. 

To summarize, for a~given input $(u(d),u(e),m)$ the enhanced machine $M_h$ computes the $m$th letter $z_m$ of $u(f)=u(d)\cdot_R\,u(e)$ 
as follows. First  $M$ disposes of the signs of $d$ and $e$ and 
determines the sign of $f$ in the obvious way. For example, if $d,e<_R 0_R$ then $d\cdot_R e=d_-\cdot_R e_-$. Then $M_h$ recovers $k''$ from $h$ 
and finds out from $h$ if 
$f\in\mathbb{R}_{\mathrm{t}}$; so we may assume that $z_m$ corresponds to a~digit $f_n$. If $f$ is terminating, $M_h$ recovers $f$ and $f_n$ from $h$ at once. Else $M$ alone determines $f_n$ as $f_n:=f(L_n)_n$, for the above defined $L_n$.
\eproof

\noindent
In the same way as in the proof of Theorem~\ref{thm_incomp} we can show that in general the order $k''$ of $d+_R e$ and $d\cdot_R e$ is incomputable.
Therefore it has to be encoded in $h$.

\section{Computing shifts}

To analyze computability of the operations $+_R$ and $\cdot_R$ in more detail, we simplify situation and  restrict these operations on $X_R$ to one variable by fixing the other at a~$d\in\R$. We call these restrictions
{\em $d$-shifts}. More generally, a~{\em unary operation} on an encoding $X$ is any map $F\colon X\to X$. {\em Computability of $F$} means computability of $F$ by a~unary machine $M$. We analyze it by means of two tools.

First, for a~unary operation $F\colon X\to X$ on an encoding $X\subset A^{\infty}$ and for $u\in X$ we say that {\em $F$ is continuous at $u$} if for every 
$k$ there is an $n_0$ such that for every $u'\in X$ coinciding with $u$ on the first $n_0$ letters, 
the value $F(u')$ coincides with $F(u)$ on the first $k$ letters. Else we say that {\em $F$ is discontinuous at $u$}. The argument in the proof of Theorem~\ref{thm_incomp}
shows in fact that if $F$ is discontinuous at some $u$ then $F$ is not computable. This agrees with the well-known result in effective 
mathematical analysis that every computable real function is continuous. 

Second, if $F\colon X\to X$ is a~unary operation on an encoding $X$, the {\em graph $D_F$ of $F$} is the directed graph
$$
D_F=(X,\,E),\ E=\{(x,\,F(x))\;|\;x\in X\}\sus X\times X\;.
$$
Its basic feature is that the only out-degree is $1$, i.e. exactly one arrow leaves every vertex.
If $F'\colon X'\to X'$ is another unary operation on another encoding $X'$, we say that {\em $F'$ is the $\pi$-conjugate of $F$} if $\pi\colon X\to X'$ is 
a~bijection and for every $u\in X'$ we have that $F'(u)=\pi(F(\pi^{-1}(u)))$. Symmetrically, $F$ is then the $\pi^{-1}$-conjugate of $F'$. This is analogous to conjugation of binary operations in the previous section. It is clear that $F$ and $F'$ are conjugates of one another if and only if they have isomorphic graphs $D_F$ and $D_{F'}$, and that conjugation maps $\pi$ are in fact isomorphisms between $D_F$ and $D_{F'}$. 

We assume that $|X|=\mathfrak{c}$ is the cardinality of continuum, as is the case for $X_R$, and mention three examples of isomorphism types of graphs $D_F$.
One is the {\em $\mathfrak{c}$-freeway} that consists of continuum many (disjoint) two-way infinite 
directed paths.Then there is the {\em $\mathfrak{c}$-sink} with  continuum many arrows heading to a~single vertex (word) which carries a~loop. Finally, {\em $\mathfrak{c}$-loops} consists of
continuum many vertices with one loop at each. 

It is easy to give examples of unary operations $F$ such that $D_F$ is the $\mathfrak{c}$-sink and $F$ is computable because $F$ is a~constant map, but not every (in fact, almost none) such $F$ is computable because $|X|=\mathfrak{c}$
but the set of (distinct types of) unary machines is countable. On the other hand, every unary operation $F$ such that $D_F$ is $\mathfrak{c}$-loops is computable because $F$ is an identity map and the unary machine computing it only needs to copy the input. It is not so obvious if there is a~unary operation $F$  such that $D_F$ 
is the $\mathfrak{c}$-freeway and $F$ is computable, but we have seen already many examples of such $F$. We characterize them in the next proposition in the case of additive $d$-shifts.

\begin{propo}\label{prop_on_add_shifts}
For any $d\in\mathbb{R}$, the additive $d$-shift
$$
u(x)\mapsto u(d+_R x),\ x\in\mathbb{R}\;,
$$
on $X_R$ is computable if and only if $d\in\mathbb{R}_{\mathrm{t}}$ and $d\le_R 0_R$.
\end{propo}
\proof
Let $d\in\mathbb{R}$ be non-terminating and $e:=1_R-_R d$. It is clear that the $d$-shift is discontinuous at $u(e)$ (see the proof of Theorem~\ref{thm_incomp}) 
and therefore is not computable. Let $d\in\mathbb{R}_{\mathrm{t}}$ with $d>_R 0_R$ and $e:=-d$. Again, the $d$-shift is discontinuous at $u(e)$ because 
the presence of the initial $-$ sign in $u(d+_R (e+_R \varepsilon))$ is sensitive to arbitrarily distant changes $\varepsilon$ of digits in $e$: for example, 
$$
1+(-1.00\ds 0100\ds)=-0.00\ds 0100\ds\;\text{ but }\;1+(-1.00\ds)=0.00\ds\;.
$$
If $d\in\mathbb{R}_{\mathrm{t}}$ with $d\le_R 0_R$ then the previous discontinuity at $u(-d)$ in sign does not occur (because of the convention (ii) in the definition of decimals).
There is no problem in computing by a~unary machine $M$ alone the sign, order and digits of $d+_R x$\,---\,see the proof of Theorem~\ref{thm_weak_comp_addi}\,---\,so 
now the $d$-shift is computable.
\eproof

Interestingly, we can make any additive $d$-shift computable by shuffling the encoding $X_R$.

\begin{thm}\label{thm_add_shifts}
For every $d\in\mathbb{R}$ there exists a~permutation (i.e., a~bijection) $$\pi\colon X_R\to X_R$$ such that the $\pi$-conjugate of the additive $d$-shift $u(x)\mapsto u(d+_R x)$ on $X_R$ is computable. 
\end{thm}
\proof
We give two proofs, the latter is more constructive than the former.

{\em A~proof by graphs. }We denote the $d$-shift as $F_d$. For $d=0_R$ the result is trivial (the graph of $F_{0_R}$ is $\mathfrak{c}$-loops). For any $d\ne 0_R$ the graph of $F_d$ is the $\mathfrak{c}$-freeway,
thus it can be conjugated by some $\pi$ to $F_{-1_R}$, which is computable by the previous proposition. A~drawback of this short proof is that 
it requires the axiom of choice: we get $\pi$ by choosing one vertex in each lane (component) of the $\mathfrak{c}$-freeways $D_{F_d}$ and $D_{F_{-1_R}}$.

{\em A~proof with explicit $\pi$. }But we can give $\pi$ explicitly. Let $d\ne0_R$. We define the permutation $\pi$ of $X_R$ by 
$$
\pi(u(e)):=u(-d^{-1}\cdot_R e)\;.
$$
We show that then the $\pi$-conjugate of the $d$-shift on $X_R$ is the computable $(-1_R)$-shift. 
We denote the $\pi$-conjugate of the $d$-shift $F_d$ on $X_R$ by $F_d'$. Then indeed
\begin{eqnarray*}
F_d'(u(e))&=&\pi(F_d(\pi^{-1}(u(e))))=\pi(F_d(u(-d\cdot_R e)))\\
&=&\pi(u(d-_R d\cdot_R e))=u(-1_R+_R e)\;.
\end{eqnarray*}
\eproof

\noindent
In view of the argument in the proof of Corollary~\ref{cor_gen_coding} it is clear that in general $\pi$ is not an equivalence of encodings.

We turn to multiplicative shifts. As for additive
shifts, we first determine those that are computable.

\begin{propo}\label{prop_on_mult_shifts}
For any $d\in\mathbb{R}$, the multiplicative $d$-shift
$$
u(x)\mapsto u(d\cdot_R x),\ x\in\mathbb{R}\;,
$$
on $X_R$ is computable if and only if $d=a_R/b_R\ge_R 0_R$ is a~nonnegative fraction such that $a\in\mathbb{N}_0$, $b\in\mathbb{N}$, $a$ and $b$ are coprime, and 
if $a\ne0$ then $a$~is not divisible by any prime number different from $2$ and $5$.
\end{propo}
\proof
For $d=0_R$ the result is trivial (the graph of the $0_R$-shift is the $\mathfrak{c}$-sink, where the sink is $u(0_R)$) 
and we assume that $d\ne0_R$. If $d<_R 0_R$ then the $d$-shift is discontinuous at $u(0_R)$ because of the sign: for example,
\begin{eqnarray*}
&&(-1)_R\cdot_R 0.00\ds 0100\ds=-0.00\ds 0100\ds\;\text{ but }\\
&&(-1)_R\cdot_R 0.00\ds=0.00\ds\;.
\end{eqnarray*}
If $d>_R 0_R$ but $d$ is not in the stated form, then it is not hard to see that $e:=d^{-1}$ is non-terminating
and, therefore, that the $d$-shift is discontinuous at $u(e)$. If $d>_R 0_R$ and $d$ is in the stated form, then the previous 
discontinuities do not occur ($d>_R 0_R$ and $d\cdot_R x\in\mathbb{R}_{\mathrm{t}}\Rightarrow x\in\mathbb{R}_{\mathrm{t}}$) and one can see, 
using the arguments in the proof of Theorem~\ref{thm_weak_comp_mult}, that a~unary machine $M$ alone computes the sign, order and digits of $d\cdot_R x$\,---\,the $d$-shift is computable.
\eproof

We consider the same problem as was solved in Theorem~\ref{thm_add_shifts} for additive shifts, if a~given multiplicative shift can be conjugated to a~computable unary operation. The answer is the same\,---\,it is always possible\,---\,but the proof is more interesting.

\begin{thm}\label{thm_mul_shifts}
For every $d\in\mathbb{R}$ there exists a~permutation $\pi\colon X_R\to X_R$ such that the $\pi$-conjugate of the multiplicative $d$-shift $u(x)\mapsto u(d\cdot_R x)$ on $X_R$ is computable. 
\end{thm}
\proof
For every $d\in\R\setminus\{-1_R,0_R,1_R\}$ the graph of the $d$-shift is the loop at $0_R$ plus the $\mathfrak{c}$-freeway, 
and since for example the $2_R$-shift is computable by the previous proposition, we see that the $d$-shifts for these decimals 
$d$ can be conjugated to a~computable unary operation, namely the $2$-shift. For $d=0_R$ and $d=1_R$ the result holds trivially
with identical map $\pi$. 

We are left with the case $d=-1_R$ when the graph $D_F$ is the loop at $0_R$ plus continuum many (disjoint) $2$-cycles. In other words,
$D_F$ is the graph of an involution on $X_R$ with exactly one fixed point at $0_R$. It suffices to find a~computable involution $F\colon X_R\to X_R$ on $X_R$ of this kind. 
We define it as follows. For the zero decimal we set $F(u(0_R))=u(0_R)$. For nonzero $d\in\R$ we set 
$F(u(d))=u(d')$, where $d'$ has the same sign as $d$ and arises from $d$ by shifting all digits $d_n$ in $d$, for all $n\in\Z$, by one position to the left (resp.
right) if the most significant nonzero digit in $d$ is at $10^n$ with odd (resp. even) $n\in\Z$. For example,
$$
d=3.14159\ds\Rightarrow d'=0.314159\ds\;\text{ and }\;d=-0.100\ds\Rightarrow d'=-1.00\ds\;.
$$
It is not hard to see that $F$ is an involution with the only fixed point $u(0_R)$ and that $F$ is computed by a~unary machine $M$. Indeed, to 
compute the $m$-th letter in $F(u(d))$ corresponding to a~digit $d_n'$ at $10^n$ in $d'$, it suffices that $M$ checks only the digits in $d$ till 
the order $10^{n-1}$.
\eproof

\noindent
The previous proof is interesting because (as we see it) at first one might think that no such computable involution $F$ exists; $F$ determines  in a way for every input decimal $d$ unambiguously if $d=0_R$ or not, in an apparent contradiction with the argument used in the proof of Theorem~\ref{thm_incomp}. But further reflection shows that here the situation is different.

\section{Concluding problems}

In the previous section we confined to two propositions and two theorems, 
but it is clear that several other intriguing questions on computability of unary and binary operations on encodings remain to be resolved. Thus
we conclude our article with two problems in this spirit. Maybe they are not too hard to solve but we do not know answers to them at the moment.

\begin{prob}
For which isomorphism types $D_F$ of graphs of unary operations $F\colon X\to X$ on encodings $X$ (with $|X|=\mathfrak{c}$) is there a~computable $F$ with the graph $D_F$?
\end{prob}

\begin{prob}
Is there a~bijection $\pi\colon X_R\to X$ such that the $\pi$-conjugate of the binary operation $+_R$ on $X_R$ is computable? The same question for $\cdot_R$.
\end{prob}

\noindent
If such $\pi$ exists then, maybe, we should seriously reconsider our decimal notation for real numbers $\dots$

\end{document}